\documentclass[12pt]{article}

\usepackage{amsmath,amssymb,amsthm}

\numberwithin{equation}{section}

\newcommand{\hb}{\hbar}
\newcommand{\pa}{\partial}
\newcommand{\Alt}{\operatornamewithlimits{Alt}}

\newtheorem{theorem}{Theorem}[section]
\newtheorem{lemma}{Lemma}[section]

\theoremstyle{definition}
\newtheorem{remark}{Remark}[section]

\newtheorem*{REMARK}{Remark}
\newtheorem*{PROPOSITION}{Proposition}

\allowdisplaybreaks

\begin{document}

\title{Quantum Geometry and Quantum Mechanics
of Integrable Systems}
\author{M.~V.~Karasev\\
\\
Moscow Institute of Electronics and Mathematics}
\date{}

\maketitle

\begin{abstract}
Quantum integrable systems and their classical counterparts are
considered. We~show that the symplectic structure and
invariant tori of the classical system can be deformed by a
quantization parameter $\hb$ to produce a new (classical)
integrable system. The~new tori selected by the $\hb$-equidistance
rule represent the spectrum of the quantum system up to~
$O(\hb^\infty)$ and are invariant under quantum dynamics in the
long-time range $O(\hb^{-\infty})$. The~quantum diffusion over the
deformed tori is described. The analytic apparatus uses quantum
action-angle coordinates explicitly constructed by an $\hb$-deformation
of the classical action-angles.
\end{abstract}

\section{Introduction}

The most elementary and fundamental systems in quantum mechanics are
integrable systems.

A quantum integrable system is a complete set of commuting operators
whose classical counterpart (arising as the quantization parameter
$\hb$ tends to zero) is a classical integrable system, i.e., admits a
complete set of functions in involution on a symplectic manifold. The joint level
surfaces of these functions are assumed to be compact, and thus they
are organized into the classical Liouville--Arnold tori fibration [1,~2].

One of basic problems is as follows: How are the spectrum and the time evolution
of a quantum integrable system related to the geometry and dynamics of the
corresponding classical integrable system?

The well-known semiclassical theory [3] enables one to compute,
approximately as $\hb\to0$, the characteristics of quantum systems
with many degrees of freedom using the objects of classical mechanics
(symplectic geometry). This geometric description works for
processes which are not
long-time and not long-distance, until the quantum
diffusion becomes significant. At long scales~${\sim\hb^{-1}}$, the
classical geometry is not adequate. In particular, in the integrable
case, the classical tori are completely destroyed under longtime
evolution by quantum diffusion. The classical symplectic
geometry fails as well in attempts to describe higher $\hb$-corrections to
the semiclassical spectrum.

As an important example, one can mention a broad class of systems
standing just in-between ``quantum'' and ``classical'' ones, like
nanosystems, where spatial sizes are comparable with the effective
wave length. For such systems, the quantization parameter $\hb$ is
not very small, say, $\hb\sim1/3$ (in dimensionless units), and thus
the critical ``long'' time or distance scale $1/\hb$ at which the
semiclassical theory fails is actually not long at all and is
practically needed to study the system. Therefore, it turns out that the classical
symplectic geometry fails to be applicable to nanosystems of this kind.

Natural questions arise. Is there a ``quantum'' symplectic geometry
which works for nanosystems? Can a quantum geometry of this kind be achieved
using an $\hb$-deformation of the classical geometry? More precisely:

{\it Is it possible to deform the original symplectic
structure and the set of functions in involution in such a way
that the corresponding deformed classical integrable system
becomes equivalent to a given quantum integrable system? Here
``equivalent'' means that the deformed tori selected by the
$\hb$-equidistance rule are invariant under the quantum dynamics and
represent the spectrum of the quantum system, as well as its
evolution.}

Actually, if we remove the word ``deformed'' from this question, then
we face a way of dealing with quantum systems as if these were classical systems
restricted to some phase-space grids subjected to quantization
rules. This coincides with Planck's basic idea [4] of the old quantum
mechanics which dominated between 1905 and 1925, before
the Heisenberg and Schr\"odinger discoveries. Thus, in a sense, the above
questions return us back to the ``naive'' concept of old quantum mechanics.

In the present paper, we show that an answer to these
questions is affirmative up to $O(\hb^\infty)$; moreover, all the
deformations mentioned above are achieved by geometrically invariant
formulas in regular domains (fibered into tori of some dimension).
Thus, we claim that the quantum and the classical
inegrability conditions are $\hb^\infty$-equivalent in regular domains. This
equivalence implies that the regular part of the spectrum, up to
$O(\hb^\infty)$, and the regular evolution at long-time range
$\sim\hb^{-\infty}$ of the quantum integrable system can be
described by using ``quantum'' geometric objects: the deformed symplectic
structure, the deformed invariant tori, the deformed frequencies of
multiperiodic rotation, and the diffusion tensor along the tori. These
objects are computed by a simple explicit algorithm.

The deformed $\hb$-equidistant tori are the actual quantum dynamical
objects replacing Planck's classical torus grids. The deformed
symplectic form vanishes on the deformed tori, and thus these
submanifolds are of Einstein--Maslov type [5, 3]. However, the
Hamiltonians of the original quantum integrable system are not
constant on these deformed tori. The actual quantum energy functions
representing the energy levels (spectrum) are obtained by a deformation
of the Hamiltonians, and the deformation takes into account the quantum
diffusion. The energy functions are in involution with respect to
the deformed Poisson bracket.

Thus, it becomes possible in a sense to replace the commutator by a
bracket to meet Dirac's ideas in quantization theory [6], as
well as to reconcile Bohr's correspondence principle with quantum
uncertainty up to $O(\hb^\infty)$.

One can conclude that, in the framework of regular integrable systems and
their perturbations, and up to accuracy $O(\hb^\infty)$, {\it quantum
mechanics turns out to be classical mechanics based on quantum
geometry}. In particular, nanosystems that are integrable or nearly
integrable can be treated up to $O(\hb^\infty)$ as classical systems with
respect to a quantum nano-geometry.

\begin{REMARK}
The term ``geometry'' is applied in this paper in
its immediate meaning, without algebraic extensions. The quantum
geometric objects are distinct from noncommutative objects or
geometric operators used in modern quantum theories. The adjective
``quantum'' is also used in its original sense, as a synonym of
``$\hb$-discretized'' and ``corresponding to the classical one as
$\hb\to0$,'' like in quantization theory [7--11] without involving
the sense of uncertainty or Hilbert space extension.
\end{REMARK}

\section{Quantum symplectic structure}

In this section, we show explicit formulas for the first
quantum corrections to classical symplectic geometry over a flat
space. The general algorithm is described in Section~4.

Assume that the algebra of quantum observables is realized by
symmetrized (Weyl ordered) functions $\widehat A=A(q,\widehat p\,)$
of the generators $q=(q^1,\dots,q^n)$ and $\widehat p=-i\hb\pa/\pa q$,
where $q^j$ stand for the Euclidean coordinates on $\Bbb R^n$. The
classical observables $A=A(q,p)$ are functions on $T^*\Bbb R^n$,
which depend on the quantization parameter in general (this
dependence is not reflected in the notation).

Let a quantum integrable system be determined by a set of commuting
self-adjoint operators $\widehat H_j$ ($j=1,\dots,n$). Thus,
\begin{equation}
[\widehat H_j,\widehat H_k]=0.
\end{equation}
In the $\hb=0$ limit, the functions
\begin{equation}
H^0_j\overset{\text{\rm def}}\to=H_j|_{\hb=0}, \qquad \{H^0_j,H^0_k\}=0,
\end{equation}
determine a classical integrable system on $T^*\Bbb R^n$ with
action-angle coordinates
$$
s=(s_1,\dots,s_n)\quad\text{and}\quad\tau=(\tau^1,\dots,\tau^n),
\qquad 0\le\tau^j\le2\pi.
$$
The classical Hamiltonians $H^0_j$ are functions on actions only.
Consider a domain in $T^*\Bbb R^n$ on which the Hamiltonians
$H^0_j$ are independent, and thus the tori
$\{s=\operatorname{const}\}$ are of the same (maximal) dimension $n$
(for details, see~[12]).

The classical symplectic form is
\begin{equation}
\omega\overset{\text{\rm def}}\to=dp\wedge dq=\frac12JdX\wedge dX, \qquad
X\overset{\text{\rm def}}\to=(\tau,s).
\end{equation}
Here $X$ is regarded as a $2n$-dimensional vector function and
$$
J=\left(
\begin{matrix}
0&I\\-I&0
\end{matrix}
\right)
$$
stands for the standard symplectic $2n\times2n$ matrix with zero
and identity $n\times n$ blocks.

Define the {\it quantum deformation of the symplectic form\/} (2.3)
up to $O(\hbar^{4})$ by the following formula:
\begin{equation}
\omega^\hb\overset{\text{\rm def}}\to=\omega+\frac{\hb^2}2J
\langle\!\langle X\overset\otimes\to{,}X\rangle\!\rangle J\,dX\wedge dX+O(\hb^4).
\end{equation}
Here the $2n\times2n$ matrix
$$\langle\!\langle X\overset\otimes\to{,}X\rangle\!\rangle$$ is composed
of the components of the vector function $X$, see~(2.3), and the
bidifferential operation
$\langle\!\langle\cdot\,,\cdot\rangle\!\rangle$ is given by
\begin{equation}
\langle\!\langle A,B\rangle\!\rangle\overset{\text{\rm def}}\to=
-\frac1{24}D^3A\cdot J\otimes J\otimes J\cdot D^3B,
\end{equation}
where $D$ stands for the derivatives with respect to the Euclidean coordinates on
$T^*\Bbb R^n$. This is the very operation entering the
$\hb$-expansion of quantum commutators,
\begin{equation}
\frac i{\hb}[\widehat A,\widehat B]
=\widehat{\{A,B\}}-\hb^2\widehat{\langle\!\langle A,B\rangle\!\rangle}+O(\hb^4).
\end{equation}

The classical form (2.3) vanishes on the classical tori
$\{s=\operatorname{const}\}$. Now we introduce new tori
annihilating the quantum form (2.4).

To each Hamiltonian $H_j$, one can assign the {\it energy function}
\begin{equation}
H^\hb_j\overset{\text{\rm def}}\to=H_j+\hb^2\Delta_sH_j+O(\hb^4) \qquad (j=1,\dots,n).
\end{equation}
Here the quantum diffusion operator $\Delta_s$ is given by
\begin{equation}
\Delta_s\overset{\text{\rm def}}\to=\frac1{16}D^2s_l\cdot J\otimes J\cdot D^2s_k\cdot
\frac{\pa^2}{\pa s_l\pa s_k}+\frac1{24}D^2s_l\cdot J\otimes J\cdot
(Ds_k\otimes Ds_m)\frac{\pa^3}{\pa s_l\pa s_k\pa s_m}
\end{equation}
(the summation over repeated Latin indices ranges from $1$ to $n$).
This is the very operator entering the general quantum composite
function expansion,
\begin{equation}
k(\widehat S)=\widehat{(I-\hb^2\Delta_S+O(\hb^4))k(S)},
\end{equation}
where the set of operators
$\widehat S=(\widehat S_1,\dots,\widehat S_n)$ on the left is
assumed to be symmetrized (Weyl ordered), see, e.g., [13,
Appendix~1, formula~(1.37)].

\begin{theorem} 
Let $\{\widehat H_j\}$ be a quantum integrable system, and let the
energy functions $H^\hb_j$ be defined by {\rm(2.7)}. Then the
deformed tori $\{H^\hb=\operatorname{const}\}$ annihilate the quantum
symplectic form $\omega^\hb$ {\rm(2.4)} up to $O(\hb^4)$. Equivalently,
the energy functions $H^\hb_j$ are in involution up to $O(\hb^4)$ with
respect to the bracket determined by the quantum Poisson tensor
$(\omega^\hb)^{-1}$.
\end{theorem}

\begin{proof}[Sketch of the proof]
Let us write out the Jacobi identity for
double commutators of the operator triple $\widehat{X^\alpha}$,
$\widehat{X^\beta}$, $\widehat{X^\gamma}$ with the vector function $X$
defined by (2.3). Taking into account formula~(2.6) and using the relation
$\{X\overset\otimes\to{,}X\}=-J=\operatorname{const}$, we obtain
\begin{equation}
\underset{\alpha,\beta,\gamma}\to{\frak S}
\{X^\alpha,\langle\!\langle X^\beta,X^\gamma\rangle\!\rangle\}=0.
\end{equation}
Here the symbol $\frak S$ stands for cyclic summation. The
Poisson bracket operation in (2.10) can be replaced by the
derivation operation in the $X$-coordinates,
$$
\underset{\alpha,\beta,\gamma}\to{\frak S}
J^{\alpha\alpha'}({\pa}/{\pa X^{\alpha'}})
\langle\!\langle X^\beta,X^\gamma\rangle\!\rangle=0,
$$
or, equivalently,
$$
\underset{\alpha,\beta,\gamma}\to{\frak S}
({\pa}/{\pa X^\alpha})J_{\beta\beta'}
\langle\!\langle X^{\beta'},X^{\gamma'}\rangle\!\rangle
J_{\gamma'\gamma}=0
$$
(the summation over repeated Greek indices ranges from $1$ to $2n$).
The last identity implies that the quantum form $\omega^\hb$ (2.4)
is closed up to $O(\hb^4)$, and thus it is symplectic.

The quantum Poisson bracket is given by
\begin{equation}
\{A,B\}^\hb\overset{\text{\rm def}}\to=dA(\omega^\hb)^{-1}\,dB
=\{A,B\}-\hb^2\frac{\pa A}{\pa X^\alpha}\langle\!\langle X^\alpha,
X^\beta\rangle\!\rangle\frac{\pa B}{\pa X^\beta}+O(\hb^4).
\end{equation}
Taking into account formula (2.1), we can derive the relation
$$
\{H_j,H_k\}=\hb^2\langle\!\langle H^0_j,H^0_k\rangle\!\rangle+O(\hb^4),
$$
and then it follows from (2.11) that
$$
\{H^\hb_j,H^\hb_k\}^\hb
=\hb^2\bigg(\langle\!\langle H^0_j,H^0_k\rangle\!\rangle
-\frac{\pa H^0_j}{\pa s_l}\langle\!\langle s_l,s_m\rangle\!\rangle
\frac{\pa H^0_k}{\pa s_m}
+\{H^0_j,\Delta_sH^0_k\}
-\{H^0_k,\Delta_sH^0_j\}\bigg)+O(\hb^4).
$$

In Lemma~2.1 below, we prove that the coefficient at $\hb^2$ on
the right-hand side vanishes, i.e.,
$$
\{H^\hb_j,H^\hb_k\}^\hb=O(\hb^4).
$$
\end{proof}

\begin{lemma} 
For the classical integrable system $\{H^0_j\}$ with the action
coordinates $\{s_j\},$ the following identities hold\/{\rm:}
\begin{equation}
\langle\!\langle H^0_j,H^0_k\rangle\!\rangle
-\frac{\pa H^0_j}{\pa s_l}
\langle\!\langle s_l,s_m\rangle\!\rangle
\frac{\pa H^0_k}{\pa s_m}
+\{H^0_j,\Delta_sH^0_k\}-\{H^0_k,\Delta_sH^0_j\}=0,
\end{equation}
where the quantum operations
$\langle\!\langle\cdot\,,\cdot\rangle\!\rangle$ and $\Delta_s$ are
determined by {\rm(2.5)} and {\rm(2.8)}.
\end{lemma}

\begin{proof}
For the first two summands in (2.12), we can see from (2.5) and
from the identities $$\{s_m,s_i\}\equiv-Ds_m\cdot J\cdot Ds_i=0$$ that
\begin{align}
&\langle\!\langle H^0_j,H^0_k\rangle\!\rangle
-\frac{\pa H^0_j}{\pa s_l}
\langle\!\langle s_l,s_m\rangle\!\rangle
\frac{\pa H^0_k}{\pa s_m}
\nonumber\\
&\qquad =\Alt_{j,k}
\bigg\{\frac{\pa^3H^0_k}{\pa s_l\pa s_m\pa s_r}
(Ds_l\otimes Ds_m\otimes Ds_r)\cdot J\otimes J\otimes J\cdot
D^3s_i\frac{\pa H^0_j}{\pa s_i}\bigg\}
\nonumber\\
&\qquad +3\Alt_{j,k}
\bigg\{\frac{\pa^2H^0_k}{\pa s_l\pa s_m}
(D^2s_l\otimes Ds_m)\cdot J\otimes J\otimes J\cdot
D^3s_i\frac{\pa H^0_j}{\pa s_i}\bigg\}
\nonumber\\
&\qquad -\frac{\pa^2H^0_j}{\pa s_l\pa s_m}
\frak S(D^2s_l\otimes Ds_m)\cdot J\otimes J\otimes J\cdot
\frak S(D^2s_r\otimes Ds_i)\frac{\pa H^0_k}{\pa s_r\pa s_i}.\quad
\end{align}
Here the symbol $\Alt$ stands for the skew-symmetric
summation
$$
\displaystyle\Alt_{j,k}\{A_{jk}\}
\overset{\text{\rm def}}\to=A_{jk}-A_{kj}
$$
and the symbol $\frak S$ for the cyclic summation
$$
{\frak S}(D^2\otimes D)_{\alpha\beta\gamma}
=D^2_{\alpha\beta}\otimes D_\gamma
+D^2_{\beta\gamma}\otimes D_\alpha
+D^2_{\gamma\alpha}\otimes D_\beta.
$$
It can readily be seen that the last summand on the right-hand side of
(2.13) vanishes.

For the third and fourth summands in (2.12), it follows from (2.8) that
\begin{align}
&\{H^0_j,\Delta_sH^0_k\}-\{H^0_k,\Delta_sH^0_j\}
\nonumber\\
&\qquad
=\Alt_{j,k}
\bigg\{\frac{\pa^3H^0_k}{\pa s_l\pa s_m\pa s_r}
D^2_{\alpha'\beta\gamma}s_l\cdot D_{\beta'}s_m\cdot D_{\gamma'}s_r\cdot
J^{\alpha'\alpha}J^{\beta'\beta}J^{\gamma'\gamma}D_\alpha s_i
\frac{\pa H^0_j}{\pa s_i}\bigg\}
\nonumber\\
&\qquad+3\Alt_{j,k}
\bigg\{
\frac{\pa^2H^0_k}{\pa s_l\pa s_m}
D^2_{\alpha'\beta'}s_l\cdot D^3_{\alpha\beta\gamma'}s_m
J^{\alpha'\alpha}J^{\beta'\beta}J^{\gamma'\gamma}D_\gamma s_i
\frac{\pa H^0_j}{\pa s_i}
\bigg\}
\nonumber\\
&\qquad +2\Alt_{j,k}
\bigg\{\frac{\pa^3H^0_k}{\pa s_l\pa s_m\pa s_r}
D^2_{\beta\gamma}s_l\cdot D^2_{\alpha'\beta'}s_r\cdot D^2_{\gamma'}s_m\cdot
J^{\alpha'\alpha}J^{\beta'\beta}J^{\gamma'\gamma}D_\alpha s_i
\frac{\pa H^0_j}{\pa s_i}
\bigg\}.
\end{align}
The second summands (with the coefficients~$3$) on the right-hand sides
of (2.13) and (2.14) are combined together into the expression
\begin{align*}
3\Alt_{j,k}
\bigg\{
&\frac{\pa^2H^0_k}{\pa s_l\pa s_m}
D^2_{\alpha'\beta'}s_l
\left(D_{\gamma'}s_mD^2_{\alpha\beta}(D_\gamma s_i)
+D^2_{\alpha\beta}(D_{\gamma'}s_m)\cdot D_\gamma s_i\right)\cdot
J^{\alpha'\alpha}J^{\beta'\beta}J^{\gamma'\gamma}
\frac{\pa H^0_j}{\pa s_i}
\bigg\}
\\
=&-3\Alt_{j,k}
\bigg\{
\frac{\pa^2H^0_k}{\pa s_l\pa s_m}
D^2_{\alpha'\beta'}s_lD^2_{\alpha\beta}\{s_m,s_i\}\cdot
J^{\alpha'\alpha}J^{\beta'\beta}
\frac{\pa H^0_j}{\pa s_i}
\bigg\}
\\
&-3\Alt_{j,k}
\bigg\{
\frac{\pa^2H^0_k}{\pa s_l\pa s_m}
D^2_{\alpha'\beta'}s_l
\left(D^2_{\alpha\gamma'}s_m\cdot D^2_{\beta\gamma}s_i
+D^2_{\beta\gamma'}s_m\cdot D^2_{\alpha\gamma}s_i\right)\cdot
J^{\alpha'\alpha}J^{\beta'\beta}J^{\gamma'\gamma}
\frac{\pa H^0_j}{\pa s_i}
\bigg\}.
\end{align*}
The last term here vanishes, as well as the bracket
$\{s_m,s_i\}=0$. Thus, the second summands in (2.13) and (2.14)
taken together give the zero value.

Thus, taking (2.13) and (2.14) together, we obtain
\begin{align}
&\langle\!\langle H^0_j,H^0_k\rangle\!\rangle
-\frac{\pa H^0_j}{\pa s_l}
\langle\!\langle s_l,s_m\rangle\!\rangle
\frac{\pa H^0_k}{\pa s_m}
+\{H^0_j,\Delta_sH^0_k\}
-\{H^0_k,\Delta_sH^0_j\}
\nonumber\\
&\qquad=\Alt_{j,k}
\bigg\{\frac{\pa^3H^0_k}{\pa s_l\pa s_m\pa s_r}
(D_{\alpha'}s_l\cdot D_{\beta'}s_m\cdot D_{\gamma'}s_r\cdot D^3_{\alpha\beta\gamma}s_i
+D^3_{\beta\gamma\alpha'}s_l\cdot D_{\beta'}s_m\cdot D_{\gamma'}s_r\cdot D_\alpha s_i
\nonumber\\
&\qquad\qquad\qquad\qquad
+2D^2_{\beta\alpha'}s_l\cdot D_{\beta'}s_m\cdot D_{\gamma'}s_r\cdot D^2_{\gamma\alpha}s_i)
J^{\alpha'\alpha}J^{\beta'\beta}J^{\gamma'\gamma}
\frac{\pa H^0_j}{\pa s_i}\bigg\}
\nonumber\\
&\qquad\qquad
-2\Alt_{j,k}
\bigg\{\frac{\pa^3H^0_k}{\pa s_l\pa s_m\pa s_r}
(D^2_{\beta\alpha'}s_l\cdot D_{\beta'}s_m\cdot D_{\gamma'}s_r\cdot D^2_{\gamma\alpha}s_i
\nonumber\\
&\qquad\qquad\qquad\qquad
-D^2_{\beta\gamma}s_l\cdot D^2_{\alpha'\beta'}s_r\cdot D_{\gamma'}s_m\cdot D_\alpha s_i)
J^{\alpha'\alpha}J^{\beta'\beta}J^{\gamma'\gamma}
\frac{\pa H^0_j}{\pa s_i}\bigg\}.
\end{align}

The last term of $\Alt_{j,k}$ on the right-hand
side of (2.15) (with the coefficient~$2$) vanishes, and
the first term of $\Alt_{j,k}$ gives the expression
$$
-\Alt_{j,k}
\bigg\{\frac{\pa^3H^0_k}{\pa s_l\pa s_m\pa s_r}
D^2_{\beta\gamma}\{s_l,s_i\}J^{\beta'\beta}J^{\gamma'\gamma}
D_{\beta'}s_m\cdot D_{\gamma'}s_r
\frac{\pa H^0_j}{\pa s_i}\bigg\}.
$$
Since $\{s_l,s_i\}=0$, we conclude that this expression, and thus the
entire expression (2.15), is identically zero.
\end{proof}

\begin{remark} 
Let us pay attention to the following important fact. The quantum
Poisson bracket (2.11), as well as the quantum symplectic form (2.4)
and the energy functions (2.7), contains not only the first derivative but also
the second and the third derivatives of phase space coordinates and Hamiltonians.
Actually, the higher-order $\hb$-corrections in (2.4), (2.7), and (2.11)
contain higher and higher derivatives (see Section~4). Objects of this kind,
depending on higher derivatives, look unusual from the viewpoint of
classical differential geometry. This is the very distinction
characterizing quantum geometry. The presence of higher
derivatives reflects the phenomena of quantum diffusion and
uncertainty. In spite of these awkward objects, the quantum geometry
is indeed a geometry, thanks to identities like (2.12), where the
quantum diffusion and dispersion operations $\Delta_s$ and
$\langle\!\langle\cdot\,,\cdot\rangle\!\rangle$ correlate with each
other and with the classical Poisson bracket.
\end{remark}

\section{Spectrum and dynamics\\ using quantum actions-angles}

By Theorem~2.1, the set of energy functions $\{H^\hb_j\}$ determines
$\operatorname{mod}O(\hb^4)$ a classical integrable system with
respect to the quantum symplectic structure $\omega^\hb$ (2.4). The
deformed tori $\{H^\hb=\operatorname{const}\}$ described in
Theorem~2.1 fiber a domain in the phase space. Hence, in this
domain, with accuracy up to $O(\hb^4)$, one
can construct the action-angle coordinates $(s^\hb,\tau^\hb)$,
$0\le\tau^{\hb j}\le2\pi$,
in the standard way. We refer to these coordinates as the {\it quantum
action-angles}. The quantum actions are constant along the deformed tori.
The explicit formula for the quantum actions is
\begin{equation}
s^\hb_j=\frac1{2\pi}\int_{C^\hb_j}\alpha^\hb, \qquad d\alpha^\hb=\omega^\hb.
\tag{3.0}
\end{equation}
Here $C^\hb_j$ is the $j$th basic cycle on the deformed torus
$\{H^\hb=\operatorname{const}\}$ containing the given phase-space point.
Let us choose an $\hb$-equidistant grid of deformed tori by taking
discrete values of the quantum actions as follows:
\begin{equation}
s^\hb=\mu+N\hb, \qquad N\in\Bbb Z^n, \quad N\sim\hb^{-1}.
\end{equation}
Here $\mu=(\mu_1,\dots,\mu_m)$ is a constant vector with the components
\begin{equation}
\mu_j=\frac{\hb}4m_j+O(\hb^2), \qquad j=1,\dots,n,
\end{equation}
where $m_j$ are the Maslov indices of the basic cycles on the tori.

\begin{theorem} 
Let $H^\hb_j$ $(j=1,\dots,n)$ be the energy functions {\rm(2.7)} of
a quantum integrable system $\{\widehat H_j\}$ {\rm(2.1),} and let
$s^\hb_j$ be the quantum actions {\rm(3.0)} determined by the
deformed integrable system $\{H^\hb_j\}$ described in
Theorem~{\rm2.1}. In this case, for a certain vector $\mu$ of type
{\rm(3.2),} the sequence of numbers
\begin{equation}
E_j[N]\overset{\text{\rm def}}\to=H^\hb_j\bigg|_{s^\hb=\mu+\hb N}, \qquad N\in\Bbb Z^n,
\end{equation}
approximates the eigenvalues of the operator $\widehat H_j$ up to
$O(\hb^4),$ i.e.,
$$
\operatorname{dist}\big(E_j[N]\,,\,\operatorname{Spectr}\widehat H_j\big)=O(\hb^4).
$$
\end{theorem}

\begin{proof}[Sketch of the proof]
The quantum symplectic form (2.4) reads
\begin{equation}
\omega^\hb=ds_j\wedge d\tau^j
+\frac{\hb^2}2\langle\!\langle s_l,s_j\rangle\!\rangle d\tau^l\wedge d\tau^j
+\frac{\hb^2}2\langle\!\langle\tau^l,\tau^j\rangle\!\rangle ds_l\wedge ds_j
+\hb^2\langle\!\langle s_j,\tau^l\rangle\!\rangle ds_l\wedge d\tau^j
+O(\hb^4).\quad
\end{equation}
On the other hand, by the definition of the quantum action-angle coordinates, we
have $\omega^\hb=ds^\hb\wedge d\tau^\hb$.
Let~us substitute into \thetag{3.4} the
$\hb$-expansions for the quantum action-angle coordinates, i.e.,
\begin{equation}
s^\hb=s+\hb^2a+O(\hb^4), \qquad
\tau^\hb=\tau+\hb^2\phi+O(\hb^4).
\end{equation}
Then it follows from (3.4) that the coefficients $a$ and $\phi$ of
expansion (3.5) satisfy the identities
$$
\frac{\pa a_j}{\pa\tau^l}-\frac{\pa a_l}{\pa\tau^j}
=\langle\!\langle s_l,s_j\rangle\!\rangle,
\qquad
\frac{\pa\phi^j}{\pa s_l}-\frac{\pa\phi^l}{\pa s_j}
=\langle\!\langle\tau^j,\tau^l\rangle\!\rangle,
\qquad
\frac{\pa a_j}{\pa s_l}+\frac{\pa\phi^l}{\pa\tau^j}
=\langle\!\langle s_j,\tau^l\rangle\!\rangle.
$$
In view of (2.6), this system implies the commutation relations
\begin{equation}
\frac i{\hb}[\widehat{s^\hb_j},\widehat{s^\hb_l}]=O(\hb^4),
\qquad
\frac i{\hb}[\widehat{\tau^{\hb\,j}},\widehat{\tau^{\hb\,l}}]=O(\hb^4),
\qquad
\frac i{\hb}[\widehat{s^\hb_j},\widehat{\tau^{\hb\,l}}]=\delta^l_j+O(\hb^4).
\end{equation}

Denote by $\operatorname{ad}_{\widehat A}$ the commutator operation
$\operatorname{ad}_{\widehat A}(\widehat B)
\overset{\text{\rm def}}\to=[\widehat A,\widehat B]$.
It follows from (3.6) that, for each~$j$, $j=1,\dots,n$, the exponential
function generated by $\frac i{\hb}\operatorname{ad}_{\widehat{s^\hb_j}}$ is
$2\pi$-periodic $\operatorname{mod}O(\hb^4)$, i.e.,
\begin{equation}
\exp\bigg\{\frac{2\pi i}{\hb}\operatorname{ad}_{\widehat{s^\hb_j}}\bigg\}=I+O(\hb^4).
\end{equation}
Since the algebra representation $A\to\widehat A$ is irreducible, it
follows from (3.7) that the exponential function generated by the
operator $\frac i{\hb}\widehat{s^\hb_j}$ is $2\pi$-periodic up to
$\operatorname{mod}O(\hb^4)$ and up to a unitary constant multiplier
\begin{equation}
\exp\Big\{\frac{2\pi i}{\hb}\widehat{s^\hb_j}\Big\}
=\exp\Big\{\frac{2\pi i}{\hb}\mu_j\Big\}\cdot I+O(\hb^4), \qquad
\mu_j=\operatorname{const}.
\end{equation}
The leading part of the Hamiltonian $s^\hb_j$ is the classical
action $s_j$ (see~(3.5)). The Hamiltonian flow generated by $s_j$
is $2\pi$-periodic. As is known [14, 15],
under this assumption, the constant $\mu_j$ in (3.8) must be of the form (3.2).

It follows from (3.8) that the spectrum of $\widehat{s^\hb_j}$
contains the sequence (3.1) (up to $O(\hb^4)$).

Note that the energy function $H_j^\hb$ depends on the quantum
action coordinates only, i.e., there exists a function $f^\hb_j$ in
$n$ variables for which
\begin{equation}
H^\hb_j=f^\hb_j(s^\hb).
\end{equation}
It follows from (2.9) and (3.9) that
$$
f^\hb_j(\widehat{s^\hb})=\widehat{H^\hb_j}
-\hb^2\widehat{\Delta_{s^\hb}H^\hb_j}+O(\hb^4).
$$
By (2.7), we have
\begin{equation}
\widehat H_j=f^\hb_j(\widehat{s^\hb})+O(\hb^4).
\end{equation}

By substituting the approximate eigenvalues (3.1) of the commuting
$\operatorname{mod}O(\hb^5)$ operators~$\widehat{s^\hb}$ into formula
(3.10), we conclude that the spectrum of $\widehat H_j$ contains the
sequence $$f^\hb_j(\mu+\hb N)+O(\hb^4).$$ In view of (3.9), this
sequence coincides with $E_j[N]+O(\hb^4)$.
\end{proof}

\begin{remark} 
The idea to compute the semiclassical spectrum using quantum
action coordinates was suggested in [16] (for details, see [9]). The
new feature in Theorem 3.1 is the explicit geometric
formula (3.0) for the quantum actions in terms of the quantum symplectic
structure and the deformed tori. Applying this formula,
one can readily compute, for instance, the {\it first quantum
corrections to classical actions\/} in expansion (3.5), namely,
\begin{equation}
a_j=\bigg(\frac{\pa H^0}{\pa s}\bigg)^{-1}
{\vphantom{\bigg(}}^l_j
\big(L_l-\langle L_l\rangle_j\big)+({1}/({2\pi}))\int_{\Sigma^0_j}\varkappa+a_{j}^{0}.
\end{equation}
Here $a_{j}^{0}$ are some constants. In (3.11) we also use the notation
$$
L_l\overset{\text{\rm def}}\to=\Delta_sH^0_l+M_l,
$$
where the symbols $M_l$ stand for the $\hb^2$-corrections in the
expansion of original Hamiltonians of the quantum integrable system,
$$
H_l=H^0_l+\hb^2M_l+\cdots\,.
$$
The angular brackets $\langle\rangle_j$ in (3.11) stand for the average
over the classical angle $0\le\tau^j\le2\pi$.
The membrane (i.e., spanning surface) $\Sigma^0_j$ in
(3.11) is bounded by the $j$th basic cycle on the classical torus
$\{H^0=\operatorname{const}\}$ containing the given phase space
point. The form $\varkappa$ in (3.11) represents the
$\hb^2$-correction to the classical form in (3.4), i.e.,
$$
\varkappa
=\frac12\langle\!\langle s_l,s_j\rangle\!\rangle d\tau^l\wedge d\tau^j
+\frac12\langle\!\langle\tau^l,\tau^j\rangle\!\rangle ds_l\wedge ds_j
+\langle\!\langle s_j,\tau^l\rangle\!\rangle ds_l\wedge d\tau^j.
$$
\end{remark}

\begin{remark} 
In view of (3.0), the {\it quantization rule\/} (3.3) can be represented
as a system of equations for the energy levels $E=E[N]$, namely,
\begin{equation}
\frac1{2\pi}\int_{C^\hb_j[E]}\alpha^\hb
=\mu_j+\hb N_j\qquad (j=1,\dots,n).
\end{equation}
Here $E=(E_1,\dots,E_n)$ are the energy levels sought for, i.e., the
eigenvalues of the commuting Hamiltonians $\widehat H_j$. The basic
cycles $C^\hb_j[E]$ belong to the deformed torus $\{H^\hb=E\}$,
where $H^\hb_j$ are the energy functions (2.7) assigned to $H_j$.
 The constants $\mu_j$ in
(3.3), (3.12) are determined using the holonomy generated by the
quantum actions following formulas (3.8). The leading term in
expansion (3.2) for the vector $\mu$ is presented by the Maslov
indices of the basic cycles on the tori, and the higher $\hb$-terms
determine certain quantum corrections.
\end{remark}

Now let us note that, by using the representation (3.10) and the
commutation relations (3.6), the quantum dynamics generated by each
Hamiltonian $\widehat H_j$ can readily be computed as follows:
\begin{equation}
\exp\Big\{({it}/{\hb})\operatorname{ad}_{\widehat H_j}\Big\}
g^{[0]}(\overset2\to{\widehat{\tau^\hb}},\overset1\to{\widehat{s^\hb}})
=g^{[t]}_j(\overset2\to{\widehat{\tau^\hb}},\overset1\to{\widehat{s^\hb}}).
\end{equation}
Here the function $g^{[0]}$ represents the initial quantum
observable or quantum density, and $g^{[t]}_j$ represents its time
evolution, i.e.,
\begin{equation}
g^{[t]}_j(\tau^\hb,s^\hb)=\exp\Big\{({it}/{\hb})
\Big(f^\hb_j\Big(s^\hb-i\hb({\pa}/{\pa\tau^\hb})\Big)
-f^\hb_j(s^\hb)\Big)+O(t\hb^4)\Big\}
g^{[0]}(\tau^\hb,s^\hb).
\end{equation}

It follows from this formula that the quantum action coordinates
are constant provided that $t\hb^4\ll1$, and all the dynamics is developed
in quantum angle coordinates only.

\begin{theorem} 
The deformed tori $\{H^\hb=\operatorname{const}\}$ are preserved by
the quantum integrable dynamics in the long-time range
$t\sim o(\hb^{-4})$. Along these tori, the dynamics generated by
$\widehat H_j$ looks like a multiperiodic rotation in the quantum
angle coordinates $\tau^\hb$ with the frequencies
$\pa H^\hb_j/\pa s^\hb$ as long as $t\sim o(\hb^{-1})$. At the time
frontier $t\sim\hb^{-1}$, the longitudinal dynamics diffuses to the
entire torus and, as long as $t\sim o(\hb^{-4}),$ this dynamics is
given by formulas {\rm(3.13)} and {\rm(3.14)}.
\end{theorem}

Note that by some transformation of the phase space near the identity,
the quantum symplectic form $\omega^\hb$ (2.4) can be turned back to
the classical form $\omega$ (2.3); however, the Hamiltonian (or the energy)
functions $H^\hb_j$ of
the deformed integrable system do not come back to $H^0_j$ under
this transformation. Thus, the Hamiltonian dynamics on deformed tori
from Theorem~3.2 is not isomorphic to the original Hamiltonian
dynamics on the classical tori.

Besides the Hamiltonian multiperiodic rotation, Theorem~3.2
describes the character of quantum spreading in the long time range.
We see that the spreading does not take place at all in directions
transversal to the deformed tori. In particular, one can claim that
there is no quantum chaos for the integrable system until
$t\sim o(\hb^{-4})$.

Along the deformed tori, the quantum spreading is explicitly
determined by formula (3.14). The leading diffusion process in (3.14)
$$
\exp\bigg\{-i\frac{t\hb}2\mathcal{D}_j\frac{\pa^2}{\pa\tau^\hb\pa\tau^\hb}\bigg\},
\qquad t\sim\hb^{-1}, \qquad
\mathcal{D}_j\overset{\text{\rm def}}\to=
\frac{\pa^2H^\hb_j}{\pa s^\hb\pa s^\hb},
$$
is easily computed via the theta function in the quantum angle
coordinates globally over each deformed torus
$$
\{s^\hb=\mu+\hb N, \quad 0\le\tau^{\hb\,l}\le2\pi \quad (l=1,\dots,n)\}.
$$

\begin{remark} 
The quantum angle coordinates are obtained geometrically from the
representation $\omega^\hb=ds^\hb\wedge d\tau^\hb$ of the quantum
symplectic structure. The explicit expression for the {\it first
quantum corrections to the classical angles} in expansion (3.5) is
the following:
\begin{equation}
\phi^l=\int^\tau_0\bigg(\langle\!\langle s_j,\tau^l\rangle\!\rangle
-\frac{\pa a_j}{\pa s_l}\bigg)\,d\tau^j+\varphi^l(s).
\end{equation}
The integration in (3.15) is taken along any path in $\tau$-space,
$a_j$ are the first quantum corrections to the classical actions
given by (3.11) and $\varphi^l$ are obtained from the initial
condition at $\tau=0$:
$$
\frac{\pa\varphi^l}{\pa s_j}-\frac{\pa\varphi^j}{\pa s_l}
=\langle\!\langle\tau^l,\tau^j\rangle\!\rangle\bigg|_{\tau=0}.
$$
The right-hand side in this condition presents coefficients of a
closed form in $s$-space (see (2.10)):
$$
\underset{k,l,j}\to{\frak S}\frac{\pa}{\pa s_k}
\langle\!\langle\tau^l,\tau^j\rangle\!\rangle=0.
$$
Thus, by the Poincar\'e lemma, one can easily resolve this initial
condition as follows
\begin{equation}
\varphi^l(s)=\int^1_0
\langle\!\langle\tau^l,\tau^j\rangle\!\rangle\bigg|_{\tau=0}
(s\xi+s^{0}(1-\xi))(s-s^{0})_{j}\xi d\xi
+\frac{\pa\psi(s)}{\pa s_l},
\end{equation}
where $s^{0}$ is some chosen point in the $s$-space and
$\psi$ is an arbitrary function in $s$-coordinates. Formula
(3.16) finally fixes the expression (3.15) for quantum corrections
to the classical angle-coordinates up to a ``gauge'' freedom
$\pa\psi(s)/\pa s_l$.

Note that (3.15) determines the correction $\phi^l$ as a
$2\pi$-periodic function in classical angles. It is guaranteed by
the identities
\begin{equation}
\big\langle\,\langle\!\langle s_j,\tau^l\rangle\!\rangle\,\big\rangle_j
=\frac{\pa}{\pa s_l}\big\langle a_j\big\rangle_j \qquad (j,l=1,\dots,n).
\end{equation}
To prove (3.17), let us note that (3.7), in the $\hb^2$-term, implies
\begin{equation}
\frac{\pa}{\pa s_l}\langle a_j\rangle_j\frac{\pa}{\pa\tau^l}
-\frac{\pa}{\pa\tau^l}\langle a_j\rangle_j\frac{\pa}{\pa s_l}
=\frac1{2\pi}\int^{2\pi}_0e^{\nu\{s_j,\cdot\}}
\langle\!\langle s_j,\cdot\rangle\!\rangle
e^{-\nu\{s_j,\cdot\}}\,d\nu.
\end{equation}
Since the left-hand side here is a first-order differential
operator, the right-hand side has to be of the same type (in spite
of the fact that $\langle\!\langle s_j,\cdot\rangle\!\rangle$ is a
third-order operator, see (2.5)). Thus, all coefficients at the
second and third derivatives on the right-hand side of (3.18)
vanish, and we derive
$$
\frac1{2\pi}\int^{2\pi}_0e^{\nu\{s_j,\cdot\}}
\langle\!\langle s_j,\cdot\rangle\!\rangle
e^{-\nu\{s_j,\cdot\}}\,d\nu
=\big\langle\,\langle\!\langle s_j,\tau^l\rangle\!\rangle\,\big\rangle_j
\frac{\pa}{\pa\tau^l}
+\big\langle\,\langle\!\langle s_j,s_l\rangle\!\rangle\,\big\rangle_j
\frac{\pa}{\pa s_l}.
$$
Comparison with the left-hand side of (3.18) implies (3.17) and so
the $2\pi$-periodicity of (3.15).
\end{remark}

\section{General algorithm for the equivalence\\
of quantum and classical integrable systems}

Let the algebra of quantum observables be realized in a space of
functions over a certain manifold~$\frak X$. The operator
representation of this algebra
\begin{equation}
A\to\widehat A
\end{equation}
is assumed to be Hermitian, invertible, and irreducible. We denote by $*$ the
product operation:
$$
\widehat{A*B}=\widehat A\widehat B,
$$
which has the properties $\overline{A*B}=\overline B*\overline A$
and $A*1=A$.

The anticommutator
$A\circledast B\overset{\text{\rm def}}\to=\frac12(A*B+B*A)$ is assumed to
admit an expansion in even powers of the quantization parameter
\begin{equation}
A\circledast B=AB-\hb^2A\odot B, \qquad
\odot=\overset{(0)}\to\odot+\hb^2\overset{(1)}\to\odot+\cdots\,.
\end{equation}
The leading term $AB$ is just the usual product of functions, and
higher $\hb$-terms in (4.2) are given by real symmetric
bidifferential operations $\overset{(\alpha)}\to\odot$ of order
$2\alpha+2$.

The commutator $[A,B]_*\overset{\text{\rm def}}\to=\frac i{\hb}(A*B-B*A)$ is
assumed to admit the $\hb$-expansion:
\begin{equation}
[A,B]_*=\{A,B\}-\hb^2\langle\!\langle A,B\rangle\!\rangle,
\qquad
\langle\!\langle\,,\,\rangle\!\rangle
=\langle\!\langle\,,\,\rangle\!\rangle^{(0)}
+\hb^2\langle\!\langle\,,\,\rangle\!\rangle^{(1)}+\cdots,
\end{equation}
whose leading term is the Poisson bracket operation
\begin{equation}
\{A,B\}=dA\,\Psi\,dB,
\end{equation}
and higher $\hb$-terms are given by real skew-symmetric
bidifferential operations
$\langle\!\langle\cdot\,,\cdot\rangle\!\rangle^{(\alpha)}$ of order
$2\alpha+3$.

Since the representation (4.1) is irreducible, the Poisson tensor
$\Psi$ (4.4) is invertible and the inverse tensor
\begin{equation}
\omega=\Psi^{-1}
\end{equation}
determines the symplectic form on the manifold $\frak X$. In
particular, the dimension of $\frak X$ is even:
$\dim\frak X=2n$.

For a set of functions $S=(S_1,\dots,S_n)$ on $\frak X$ and a
function $k$ of $n$ variables, one can define the symmetric (Weyl
ordered) $*$-composite function $k(S)_*$, as well as the usual
composite function $k(S)$. Their difference can be expanded into an
$\hb$-power series to $O(\hb^\infty)$:
\begin{equation}
k(S)_*=(I-\hb^2\Delta_S)k(S), \qquad
\Delta_S=\Delta_S^{(0)}+\hb^2\Delta_S^{(1)}+\cdots\,.
\end{equation}
The coefficients of expansion (4.6) are explicitly computed via the
operations $\overset{(\alpha)}\to\odot$ from (4.2). For instance, the
leading ``diffusion'' term in (4.6) reads
\begin{equation}
\Delta_S^{(0)}=\frac12S_j\overset{(0)}\to\odot S_l\frac{\pa^2}{\pa S_j\pa S_l}
+\frac16[[S_j\overset{(0)}\to\odot,S_l],S_m]\frac{\pa^3}{\pa S_j\pa S_l\pa S_m}.
\end{equation}

Let $H_1,\dots,H_n$ be the Hamiltonians of a {\it quantum integrable
system}. Operators $\widehat H_j$ mutually commute, and so the
functions $H_j$ commute with respect to the $*$-product:
\begin{equation}
[H_j,H_k]_*=0.
\end{equation}

Let this set of functions admit the $\hb$-expansion to
$O(\hb^\infty)$:
\begin{equation}
H_j=H^0_j+\hb^2M_j, \qquad
M_j=M^{(0)}_j+\hb^2M^{(1)}_j+\cdots\,.
\end{equation}

The leading term of (4.9) determines the classical integrable
system, i.e., the functions in involution with respect to the
Poisson bracket (4.4). We assume that the level surfaces of these
functions are compact.
In a regular domain, where all $H_j^0$ are independent, the dimension
of these surfaces (tori) is $n=\frac12\dim\frak X$.

Denote by $s_1,\dots,s_n$ classical action coordinates on this
regular domain, and let $\tau^1,\dots,\tau^n$ be the corresponding
angle coordinated so that $H^0_j=f_j(s)$ and
\begin{equation}
\{s_j,s_l\}=0, \qquad
\{\tau^j,\tau^l\}=0, \qquad
\{s_j,\tau^l\}=\delta^l_j.
\end{equation}

In order to set up the equivalence of the quantum integrable system
$\{\widehat H_j\}$ to a classical (deformed) integrable system, we
first introduce the main notation and then explain how the
algorithm of equivalence works.

The quantum action-angle coordinates $s^\hb$, $\tau^\hb$ are
determined by $\hb$-expansions to $O(\hb^\infty)$
\begin{equation}
\alignedat2
s^\hb&=s+\hb^2a, &\qquad a&=a^{(0)}+\hb^2a^{(1)}+\cdots,
\\
\tau^\hb&=\tau+\hb^2\phi, &\qquad \phi&=\phi^{(0)}+\hb^2\phi^{(1)}+\cdots,
\endalignedat
\end{equation}
to obey the commutation relations
\begin{equation}
[s^\hb_j,s^\hb_l]_*=O(\hb^\infty),
\qquad
[\tau^{\hb\,j},\tau^{\hb\,l}]_*=O(\hb^\infty),
\qquad
[s^\hb_j,\tau^{\hb\,l}]_*=\delta^l_j+O(\hb^\infty).
\end{equation}
The choice of quantum actions must be consistent with the given
quantum integrable system (4.8). Namely, the Hamiltonians $H_j$ must
be expressed as $*$-composite functions in quantum actions to
$O(\hb^\infty)$:
\begin{equation}
H_j=f^\hb_j(s^\hb)_*\,, \qquad
f^\hb_j=f_j+\hb^2g_j, \qquad
g_j=g^{(0)}_j+\hb^2g^{(1)}_j+\cdots\,.
\end{equation}

We associate each Hamiltonian $H_j$ with the energy function
\begin{equation}
H^\hb_j\overset{\text{\rm def}}\to=f^\hb_j(s^\hb), \qquad
H^\hb_j=H^0_j+\hb^2L_j, \qquad
L_j=L^{(0)}_j+\hb^2L^{(1)}_j+\cdots\,.
\end{equation}
In view of formula (4.6), we have the following relation:
$$H_j=(I-\hb^2\Delta_{s^\hb})H^\hb_j,$$ and hence
\begin{equation}
H^\hb_j=(I-\hb^2\Delta_{s^\hb})^{-1}H_j
=\Big(I+\sum_{k\ge1}\hb^{2k}(\Delta_{s^\hb})^k\Big)H_j.
\end{equation}
And finally, by means of the quantum action-angle coordinates, we
define the quantum symplectic form
\begin{equation}
\omega^\hb=ds^\hb\wedge d\tau^\hb, \qquad
\omega^\hb=\omega+\hb^2\varkappa, \qquad
\varkappa=\varkappa^{(0)}+\hb^2\varkappa^{(1)}+\cdots,
\end{equation}
as well as the quantum Poisson tensor $\Psi^\hb=(\omega^\hb)^{-1}$
and the corresponding quantum bracket $\{A,B\}^\hb=dA\Psi^\hb dB$.

In order to construct all these quantum geometric objects, we
proceed by induction.

At the zero induction step, the coefficients $a^{(0)}$,
$\phi^{(0)}$, $g^{(0)}_j$, $L^{(0)}_j$, $\varkappa^{(0)}$ of
$\hb$-expansions (4.11), (4.13), (4.14), (4.16) are known from the
procedure described in Section~2. Assume that we know all
coefficients $a^{(\alpha)}$, $\phi^{(\alpha)}$, $g^{(\alpha)}_j$,
$L^{(\alpha)}_j$, $\varkappa^{(\alpha)}$ for $\alpha\le k-1$, and
let us demonstrate how to compute them for $\alpha=k$.

The coefficients $L^{(k)}_j$ are computed by (4.15) using expansions
(4.9). Thus the deformed tori $\{H^\hb=\operatorname{const}\}$
are known at the order $\hb^{2k}$.

The form $\varkappa^{(k)}$ in (4.16) is given by
\begin{align}
\hskip-1cm
\varkappa^{(k)}
=\frac12\bigg(\frac{\pa a^{(k)}_j}{\pa\tau^l}
&-\frac{\pa a^{(k)}_l}{\pa\tau^j}\bigg)
\,d\tau^l\wedge d\tau^j
-\frac12\bigg(\frac{\pa\phi^{(k)j}}{\pa s_l}-\frac{\pa\phi^{(k)l}}{\pa s_j}\bigg)
\,ds^l\wedge ds^j
\nonumber\\
&+\bigg(\frac{\pa a^{(k)}_l}{\pa s_j}+\frac{\pa\phi^{(k)j}}{\pa\tau^l}\bigg)
\,ds_j\wedge d\tau^l
+\!\!\!\!\sum_{0\le m\le k-1}\!\! da^{(m)}_j\wedge d\phi^{(k-m-1)j}.
\end{align}
The combinations of first derivatives of functions $a^{(k)}$ and
$\phi^{(k)}$ standing in first three summands are derived from
commutation relations (4.12) at the order $\hb^{2k}$. For instance,
\begin{align*}
\frac{\pa a^{(k)}_j}{\pa\tau^l}-\frac{\pa a^{(k)}_l}{\pa\tau^j}
&=\langle\!\langle s_l,s_j\rangle\!\rangle^{(k)}
+\!\!\!\!\sum_{0\le m\le k-1}\{a^{(m)}_j,a^{(k-m-1)}_l\}
\\
&\qquad +\!\!\!\!\sum_{0\le m\le k-1}
\left(\langle\!\langle s_l,a^{(k-m-1)}_j\rangle\!\rangle^{(m)}
-\langle\!\langle s_j,a^{(k-m-1)}_l\rangle\!\rangle^{(m)}\right)
\\
&\qquad +\!\!\!\!\sum_{0\le m\le k-2}\,\,\sum_{0\le r\le k-m-2}
\!\!\!\! \langle\!\langle a^{(r)}_l,a^{(k-m-r-2)}_j\rangle\!\rangle^{(m)}.
\end{align*}
Analogous equations hold for
${\pa\phi^{(k)j}}/{\pa s_l}-{\pa\phi^{(k)l}}/{\pa s_j}$ and
for
${\pa a^{(k)}_l}/{\pa s_j}+{\pa\phi^{(k)j}}/{\pa\tau^l}$.
The right-hand sides of these equations are known by the inductive
hypothesis. By substituting these right-hand sides into the first
three summands in (4.17), we derive a formula for
$\varkappa^{(k)}$ via $a^{(\alpha)}$ and $\phi^{(\alpha)}$ with
$\alpha\le k-1$. Thus we compute the quantum form $\omega^\hb$
(4.16) at the order $\hb^{2k}$.

Then it has to be checked that the form $\omega^\hb$ vanishes on the
deformed tori up to $O(\hb^{2k+2})$.

As soon as we know a symplectic structure and a tori fibration
annihilating the symplectic form up to $O(\hb^{2k+2})$, the
action-angle coordinates are automatically constructed with the same
accuracy. Therefore, we have obtained $a^{(k)}$ and $\phi^{(k)}$ as well as the
coefficient $g^{(k)}$ via expansions (4.14), (4.13). The induction
is completed.

\begin{PROPOSITION}
Theorems~{\rm 2.1, 3.1,} and~{\rm3.2} hold with accuracy
$O(\hb^\infty)$ {\rm(}instead of $O(\hb^4)${\rm)} and at the time
range $t\sim O(\hb^{-\infty})$ {\rm(}instead of $O(\hb^{-4})${\rm)}
in the general setting described in the present section.
\end{PROPOSITION}

It is important to stress that the algorithm for
$\hb^\infty$-equivalence of quantum and classical integrability
described at this section works in arbitrary quantization schemes
over general symplectic manifolds. Thus this algorithm can be
applied to a wide variety of integrable systems appearing in
different areas of mathematical physics (see 17--22]).

\end{document}